\documentclass[12pt]{article}
\usepackage{amssymb,amsfonts,amsmath, psfrag,eepic,colordvi,graphicx,epsfig}
\usepackage{amssymb,latexsym,graphics,array}
\usepackage{MnSymbol}
\usepackage[enableskew]{youngtab}
\usepackage{float,enumerate,enumitem}
\usepackage{tikz,ifthen}
\usepackage{pgfplots}
\usepackage{hyperref}
\usepackage{booktabs}
\usepackage{mathtools}
\usetikzlibrary{math}

\hypersetup{colorlinks=true}%give color to the clickable link
\allowdisplaybreaks[4]%%allow display breaks

\topskip  
\newcommand{\rmnum}[1]{\romannumeral #1}

\numberwithin{equation}{section}
\newtheorem{theorem}{Theorem}[section]
\newtheorem{proposition}[theorem]{Proposition}

\newtheorem{lemma}[theorem]{Lemma}

\begin{document}
	\parskip 6pt
	
	\pagenumbering{arabic}
	\def\sof{\hfill\rule{2mm}{2mm}}
	\def\ls{\leq}
	\def\gs{\geq}
	\def\SS{\mathcal S}
	\def\qq{{\bold q}}
	\def\MM{\mathcal M}
	\def\TT{\mathcal T}
	\def\EE{\mathcal E}
	\def\lsp{\mbox{lsp}}
	\def\rsp{\mbox{rsp}}
	\def\pf{\noindent {\it Proof.} }
	\def\mp{\mbox{pyramid}}
	\def\mb{\mbox{block}}
	\def\mc{\mbox{cross}}
	\def\qed{\hfill \rule{4pt}{7pt}}
	\def\block{\hfill \rule{5pt}{5pt}}
	\def\Asctop{\mathrm{Asctop}}
	\def\Ascbot{\mathrm{Ascbot}}
	\def\Desbot{\mathrm{Desbot}}
	\def\Destop{\mathrm{Destop}}
	\def\asctop{\mathrm{asctop}}
	\def\ascbot{\mathrm{ascbot}}
	\def\desbot{\mathrm{desbot}}
	\def\destop{\mathrm{destop}}
	\def\A{\mathcal{A}}
	\def\MA{\hat{\mathcal{A}}}
	\def\RA{\hat{\mathcal{B}}}
	\def\CC{\hat{\mathcal{C}}}
	\def\DD{\hat{\mathcal{D}}}
	\def\nub{\text{Nub}}
	\def\Cay{\text{Cay}}
	\def\hh{\text{hat}}
\def\ha{\eta}
	\def\Add{\mathrm{Add}}
	\def\std{\mathrm{std}}
	\def\T{\mathcal{T}}

	\def\lr#1{\multicolumn{1}{|@{\hspace{.6ex}}c@{\hspace{.6ex}}|}{\raisebox{-.3ex}{$#1$}}}
	\def\red{\textcolor{red}}

	\begin{center}
{\Large \bf Pattern avoidance in revised ascent sequences}
	\end{center}
	
	\begin{center}
		{\small  Robin D.P. Zhou}

		College of Mathematics Physics and Information\\
		Shaoxing University\\
		Shaoxing 312000, P.R. China\\

Email: dapao2012@163.com

%		$^{b}$Department of Mathematics,
%		Zhejiang Normal University\\
%		Jinhua 321004, P.R. China		
	\end{center}
	
\noindent {\bf Abstract.} 	
	Inspired by the definition of modified ascent sequences, we introduce a new class of integer sequences called revised ascent sequences. These sequences are defined as Cayley permutations where each entry is a leftmost occurrence if and only if it serves as an ascent bottom.
	We construct a bijection between ascent sequences and revised ascent sequences by adapting the classic hat map, which transforms ascent sequences into modified ascent sequences. 
	Additionally, we investigate revised ascent sequences that avoid a single pattern, leading to a wealth of enumerative results. Our main techniques include the use of bijections, generating trees, generating functions, and the kernel method.

	\noindent {\bf Keywords}: ascent sequence, revised ascent sequence, pattern avoidance, generating tree,  kernel method.

	\noindent {\bf AMS  Subject Classifications}: 05A05, 05C30

	%===========================================================================
	
	\section{Introduction}
	
%Let  $x = x_1x_2\cdots x_n$ be a sequence of positive integers.
%%Given a sequence  $x = x_1x_2\cdots x_n$ of positive integers and an 
%%index $i$,  if $x_i > x_{i+1}$ 
%An index $i$ ($1\leq i \leq n-1$) is said to be an \emph{ascent bottom}
%(resp. \emph{descent top}) of $x$ if $x_{i}<x_{i+1}$ (resp. $x_{i}>x_{i+1}$).
%Similarly, an index $i$ ($2\leq i \leq n$) is said to be an \emph{ascent top}
%(resp. \emph{descent bottom}) of $x$ if $x_{i-1}<x_{i}$ (resp. $x_{i-1}>x_{i}$).
%We always designate $1$ as an ascent top (resp. ascent bottom, descent top, descent bottom) of $x$ if $x$ is non-empty.

Let $n$ be a positive integer and let $[n] = \{1,2,\ldots,n\}$.
An {\em endofunction} of size $n$ is a map $x:[n]\mapsto [n]$, which 
is represented  as the word $x = x_1x_2\cdots x_n$ with $x_i = x(i)$.
For an endofunction $x = x_1x_2\cdots x_n$, let 
\begin{align*}
	\Asctop(x) &:=  \{1\} \cup \{i \mid 2\leq i\leq n \text{ and } x_{i-1} < x_i \},\\
	\mathrm{Ascbot}(x) &:=   \{1\} \cup\{i \mid 1\leq i\leq n-1 \text{ and } x_{i} < x_{i+1} \},\\
		\mathrm{Destop}(x) &:=   \{1\} \cup\{i \mid 1\leq i\leq n-1 \text{ and } x_{i} > x_{i+1} \} \quad \text{and}\\
	\mathrm{Desbot}(x) &:=  \{1\} \cup \{i \mid 2\leq i\leq n \text{ and } x_{i-1} > x_{i} \}
\end{align*}
be the sets of ascent tops, ascent bottoms, descent tops and descent bottoms of $x$, respectively.
Depending on the context, these definitions can also refer to the values rather than the indices if no confusion arises.
For example, if $x_i< x_{i+1}$, we also call $x_i$ an ascent bottom 
and $x_{i+1}$ an ascent top.
Note that we always include $1$ as an element in each of these four sets. 
Let $\asctop(x)$ (resp. $\ascbot(x)$, $\destop(x)$ and $\desbot(x)$)
be the cardinality of $\Asctop(x)$ (resp. $\Ascbot(x)$, $\Destop(x)$ and $\Desbot(x)$).
%\[\asctop(x) = |\Asctop(x)|, \ascbot(x) = |\Ascbot(x)|\]
%$asctop(x)$ 
%This assumption simplifies the definitions of modified ascent sequences and revised ascent sequences

We call $x$ an \emph{ascent sequence} if $x_1=1$ and $1 \leq x_i \leq \asctop(x_1x_2\cdots x_{i-1})+1$ for all $2 \leq i \leq n$.
An example of such a sequence is $12213245$, whereas $112142$ is not, because the $4$ is greater than $\asctop(1121) +1= 3$.
Let $\mathcal{A}_n$ denote the set of ascent sequences of length $n$.
Ascent sequences were introduced by Bousquet-M\'elou, Claesson, Dukes and Kitaev \cite{Bousquet} to build connections among 
three other seemingly unrelated  combinatorial objects:
$(2+2)$-free posets, a family of
permutations avoiding a certain pattern and Stoimenow matchings \cite{Stoimenow}.
All these combinatorial objects were proved to be enumerated by the 
\emph{Fishburn numbers}, appearing as sequence A022493 in the OEIS \cite{oeis} with the following elegant generating function:
\begin{equation*}
	\sum_{n\geq1}\prod_{i=1}^n(1-(1-x)^i)=x+2x^2+5x^3+15x^4+53x^5+217x^6+1014x^7+\cdots.
\end{equation*}

The combinatorial objects which are enumerated by the Fishburn numbers
are called \emph{Fishburn structures}.
Ascent sequences have since been studied in a series of papers, 
connecting them to many newly discovered Fishburn structures, such as 
Fishburn matrices \cite{Yan2019,Dukes2010}, descent correcting 
sequences~\cite{Claesson2},  Fishburn diagrams~\cite{Levande}, inversion sequences avoiding a vincular pattern of length $3$~\cite{LY}, Fishburn covers and Fishburn trees \cite{Cerbai23a}.
The study of ascent sequences have been generalized to weak ascent
sequences \cite{Benyi,Benyi2024} and 
$d$-ascent sequences \cite{Dukes2023,Zhou}.
These connections and generalizations have established deep connections in these apparently disparate combinatorial structures, including 
permutations, certain integer matrices, posets, trees, and set partitions, see also \cite{Chen2013,Dukes2019,Lin,Yan} and references therein.

\emph{Modified ascent sequences} are originally defined as 
the bijective images of ascent sequences under a certain hat map \cite{Bousquet}.
They, along with ascent sequences, have played a crucial role in the study of Fishburn structures.
A \emph{Cayley permutation} is an endofunction $x$ such that 
$x$ contains each integer between $1$ and $\max(x)$.
Let $\Cay_n$ denote the set of Cayley permutations of length $n$.
Very recently, 
Cerbai and Claesson~\cite{Cerbai23a} provided
an alternative definition of modified ascent sequences as 
Cayley permutations where each entry is a leftmost copy if and only if it is an ascent top.
Let 
$\nub(x)$ denote the set of indices of the leftmost copies of $x$. 
That is, 
\[\nub(x) = \{\min \,x^{-1}(i) \mid  i \in \mathrm{Im}(x)\},\]
where $\mathrm{Im}(x)$ is the set of images of $x$.
Then the set $\MA_n$ of modified ascent sequences of length $n$ is defined as
\[\MA_n = \{x\in \text{Cay}_n \mid \Asctop(x) = \nub(x)\}.\]
For example, $x = 135144312$ is a modified ascent sequence because
$\Asctop(x) = \nub(x) = \{1,2,3,5,9\}$.
Subsequently, modified ascent sequences have been extensively studied by Cerbai, Claesson, and Sagan \cite{Cerbai23b,Cerb,Cera,Cerb2406,Cerb2408}.

Inspired by the definition of modified ascent sequences, we are 
motivated to investigate the following sets:
\begin{align*}
	\RA_n &:=  \{x\in \text{Cay}_n \mid \Ascbot(x) = \nub(x)\},\\
	\CC_n &:=  \{x\in \text{Cay}_n \mid \Destop(x) = \nub(x)\},\\
	\DD_n &:=  \{x\in \text{Cay}_n \mid \Desbot(x) = \nub(x)\}.
\end{align*}
We will demonstrate that all these sets are enumerated by the Fishburn numbers.
Given an endofunction $x = x_1x_2\cdots x_n$ with $\max(x) = m$, its \emph{complement},
denoted by $x^c$, is defined by $x^c(i) = m+1 - x_i$ for $1\leq i\leq n$.
It is clear that $\nub(x) = \nub(x^c)$, $\Asctop(x) = \Desbot(x^c)$, and $\Ascbot(x) = \Destop(x^c) $.
Therefore, the complement mapping induces a bijection between 
$\MA_n$ and $\DD_n$, as well as a bijection between $\RA_n$ and $\CC_n$.
It follows that $|\RA_n| = |\CC_n|$ and $|\DD_n| = |\MA_n| = F_n$, where $F_n$ denotes the $n$-th Fishburn number.

The objective of this paper is to study the set $\RA_n$, whose elements we refer to as \emph{revised ascent sequences}. 
We have chosen this name because these sequences exhibit behavior similar to that of modified ascent sequences, which can be obtained from ascent sequences through a specific hat map. 
Additionally, this naming convention serves to distinguish them from traditional modified ascent sequences.
After establishing a bijection between ascent sequences and 
revised ascent sequences in Section \ref{sec:bijection},
we will focus on studying revised ascent sequences that avoid a single pattern.
This work parallels the research  of Duncan and Steingr\'imsson  on ascent sequences~\cite{Duncan} and  Cerbai's work on modified 
ascent sequences \cite{Cera}.

Let $x\in \Cay_n$ and $\sigma \in \Cay_k$ be two Cayley permutations 
with $k\leq n$.
An \emph{occurrence} of $\sigma$ in $x$ is a subsequence $x_{i_1}x_{i_2}\cdots x_{i_k}$ which is order isomorphic 
to $\sigma$. 
Here, \emph{order isomorphic} means that $x_{i_s}<x_{i_t}$ if and only if
$\sigma_s < \sigma_t$ and $x_{i_s}=x_{i_t}$ if and only 
if $\sigma_s = \sigma_t$.
We say $x$ \emph{contains} the pattern $\sigma$ if $x$ contains an occurrence of $\sigma$.
Otherwise, we say that $x$ \emph{avoids} the pattern $\sigma$ or $x$ is \emph{$\sigma$-avoiding}.
For example,  the revised ascent sequence $4134232$ contains two
occurrences of the pattern $123$, namely, the subsequences $134$ and $123$, but it avoids the pattern $112$.
Let $\RA_n(\sigma)$ denote the set of revised ascent sequences of length $n$ that avoid the pattern $\sigma$ and let $$\RA_{\sigma}(t)= \sum_{n\geq 1}|\RA_n(\sigma)|t^n$$
be the generating function of $\RA_n(\sigma)$.

 The rest of this paper is organized as follows.

In Section \ref{sec:bijection}, we prove that revised ascent sequences are enumerated by 
the Fishburn numbers by constructing a bijection between 
$\A_n$ and $\RA_{n+1}$. 
%Our bijection is analogous to the classic hat map from ascent sequences
%to modified ascent sequences, with the 
%difference being that, when go through each ascent,  we increment by one all the entries in the corresponding prefix that are currently greater than or equal to the current ascent bottom, rather than the ascent top.
We also provide a generating tree for revised ascent sequences,
which shares the same succession rules with the generating tree for 
ascent sequences.

In Section \ref{sec:pattern},  we investigate the enumeration of revised ascent sequences that avoid a single pattern.  
%using various methods, including bijections, generating trees,  generating functions and the kernel method.
Our research has produced a rich set of enumerative results.
An overview of these findings can be found in Table \ref{table:main-results}, 
where the patterns in the same row determine the same set of sequences.
%The question mark indicates that the enumeration of this pattern has not yet been determined.

\begin{table}[t]  
	\centering  
	\begin{tabular}{>{\centering\arraybackslash}p{4cm} >{\centering\arraybackslash}p{5cm}  >{\centering\arraybackslash}p{4cm}}  
		\toprule  
		$\sigma$ & ${|\RA_n(\sigma)|}_{n\geq 2}$ & Reference \\
		\midrule  
		11 &  0, 0, 0,$\ldots$ & Section \ref{subsection:12} \\  
		12, 21, 212 & 1, 1, 1,$\ldots$ & Section \ref{subsection:12} \\  
		221&  1, 2, 3, 4,$\ldots$ & Section \ref{subsection:221} \\
		312, 122 & $2^{n-2}$ &  Sections \ref{subsection:312}, \ref{subsection:122} \\    
		231, 321, 3231 & $2^{n-1}-n+1$ &  Section \ref{subsection:321} \\  
		213 & Catalan number &   Section \ref{subsection:213}\\
		121, 211, 2121 & A026898 &  Section \ref{subsection:211}\\
		112 & Bell number &   Section \ref{subsection:112}\\
		123 & A105633 &  Section \ref{subsection:123} \\  
		132, 3132 & A082582 &   Section \ref{subsection:132}\\
		111 & 1, 1, 2, 4, 10, 29, 97,\ldots &   open problem\\
		\bottomrule  
	\end{tabular}  
	\\
	\caption{Enumeration of revised ascent sequences avoiding a single pattern $\sigma$.} \label{table:main-results}
\end{table}

Section~\ref{sec:final} contains some final remarks and suggestions for future work.

\section{A bijection between $\A_n$ and $\RA_{n+1}$}\label{sec:bijection}

% with the 
%difference being that, when go through each ascent,  we increment by one all the entries in the corresponding prefix that are currently greater than or equal to the current ascent bottom, rather than the ascent top.
%\begin{enumerate}
%	\item Do for $i=a_1,a_2,\ldots,a_k$ \\\hspace*{0.5cm}
%	      \hspace*{1cm}  for $j=1,2,\ldots, i-1$ \\
%	       \hspace*{1.5cm} if $x_j\geq  x_{i}$ then $x_j:=x_j+1$.
%	       Let $\tilde{x}$ be the resulting sequence.
%	 \item  Add the maximum value of $\tilde{x}$ at the beginning of 
%	  $\tilde{x}$.
%\end{enumerate}

This section is devoted to constructing a bijection $\ha$ between 
ascent sequences and revised ascent sequences.
Our bijection closely mirrors the classic hat map that transforms ascent sequences into modified ascent sequences,
but with a notable difference in the modification process.

Let $x$ be an ascent sequence with $\Ascbot(x) = \{a_1,a_2,\ldots, a_k\}$, where $a_1<a_2<\cdots <a_k$.
Then $\ha(x)$ is defined by the following procedure: 
%\begin{center}
%		 for $i=a_1,a_2,\ldots,a_k$ do\\
%	\hspace*{1cm}  for $j=1,2,\ldots, i-1$ do\\
%	\hspace*{3cm} if $x_j\geq  x_{i}$ then $x_j:=x_j+1$.
%\end{center}
{\par
\centering
\hspace*{-1.5cm}for $i=a_1,a_2,\ldots,a_k$ do\\
  for $j=1,2,\ldots, i-1$ do\\
\hspace*{2cm} if $x_j\geq  x_{i}$ then $x_j:=x_j+1$.\par
}
\noindent Define $\tilde{x}$ to be the resulting sequence after performing the above operations on $x$.
Then $\ha(x)$ is the sequence obtained from $\tilde{x}$ by 
prepending the maximum value of $\tilde{x}$ to the beginning.
  Observe that if we replace $\Ascbot(x)$ with $\Asctop(x)$ in the above construction, 
   $\tilde{x}$ becomes the corresponding modified ascent sequence
  of $x$ under the classic hat map.
As an example, consider the ascent sequence $x = 12132124$.
We have $\Ascbot(x) = \{1,3,6,7\}$ and $\tilde{x}$ is obtained through the following steps:
\begin{align*}
	x =& \red{1}\;2 \;1 \;3 \;2 \;1\; 2 \;4 \\
  	  & 1\; 2 \;\red{1} \;3 \;2 \;1\; 2 \;4 \\
  	  &  2\; 3 \;1 \;3 \;2 \;\red{1}\; 2 \;4 \\
  	  &  3\; 4 \;2 \;4 \;3 \;1\; \red{2} \;4 \\
  	  &  4\; 5 \;3 \;5 \;4 \;1\; 2\;4 = \tilde{x}.
\end{align*}
In each step, every entry strictly to the left of and weakly larger than the number colored by red is incremented by one.
Then we have $\ha(x) = 545354124$.
Observe that $\Ascbot(x) = \Ascbot(\tilde{x})$.
The construction of the map $\ha$ can be easily inverted.
Hence, the map $\ha$ is injective.

Let $x = x_1x_2\cdots x_n$ be an endofunction and let $1\leq \ell\leq \max(x) +1$.
Define $y = \Add(x,\ell)$ to be the following sequence of length $n+1$ depending on 
whether $x_n<\ell$.
\begin{itemize}
	\item If $1\leq \ell \leq x_n$, we have $y = x\ell$. Here $x\ell$ is the concatenation of $x$ and $\ell$.
	\item If $x_n <\ell\leq    \max(x)+1$,  we have $y_n = x_n$, $y_{n+1} = \ell$ and for $1\leq i \leq n-1$,
	\[y_i=
	\begin{cases}
		x_i,     \quad &\text{if }  x_i < x_n, \\
		x_i +1,  \quad&\text{if }  x_i \geq x_n.
	\end{cases}\]
\end{itemize}
It is obvious that $\tilde{x} = \Add(\tilde{x'},x_n)$, where 
$x' = x_1x_2\cdots x_{n-1}$.

Conversely, one can easily recover $x$ from $y$.
Let $y = y_1y_2\cdots y_{n+1}$ be an endofunction.
Define $x = \text{Rom}(y)$ to be the following sequence of length $n$ depending on whether 
$y_n <y_{n+1}$.
\begin{itemize}
	\item If $1\leq y_{n+1} \leq y_n$, we have $x = y_1y_2\cdots y_n$. 
	\item If $y_n <y_{n+1}$,  we have $x_n = y_n$, and for $1\leq i \leq n-1$,
	\[x_i=
	\begin{cases}
		y_i,     \quad &\text{if }  y_i \leq  y_n, \\
		y_i -1,  \quad&\text{if }  y_i > y_n.
	\end{cases}\]
\end{itemize}
Then we have $\text{Rom}(\Add(x,\ell)) = x$.
Additionally, there exists a relation between $\ha$ and $\Add$ as follows.
%From the construction of the map $\ha$, we have 
%%$$\tilde{x} = \Add(x_1,x_2,\ldots,x_n),$$
%$\tilde{x} = \Add(\tilde{x'},x_n)$,
% where $x' = x_1x_2\cdots x_{n-1}$.
%$$\Add(x1,x2,\ldots,x_n) = \Add(\Add(x_1,x_2,\ldots,x_{n-1}),x_n).$$

\begin{lemma}\label{lem:Add}
	Let $n\geq 2$ and let $x = x_1x_2\cdots x_n$ be an ascent sequence.
	Then we have $\asctop(x) = \max(\ha(x))$ and $\ha(x) = \Add(\ha(x'),x_n)$,  where $x' = x_1x_2\cdots x_{n-1}$.
%	$$\ha(x) = \Add(1,x_1,x_2,\ldots,x_n).$$
\end{lemma}

\pf
Let us first focus on the property $\asctop(x) = \max(\ha(x))$.
By the construction of the map $\ha$, we have $ \max(\ha(x)) = \max(\tilde{x})$.
Thus, to prove $\asctop(x) = \max(\ha(x))$, it is sufficient to prove 
$\asctop(x) = \max(\tilde{x})$.
We will prove this result by induction on $n$ where $n=2$ is trivial.
Assume the result for $n-1$. 
We need to prove the result for $n$.
%Let $x' = x_1x_2\cdots x_{n-1}$. 
By the induction hypothesis, we have $\asctop(x') = \max(\tilde{x'})$.
Now we consider two cases.
If $1\leq x_n \leq x_{n-1}$, 
then $\tilde{x} =\Add(\tilde{x'},x_n) = \tilde{x'}x_n$.
Note that  $x_n\leq x_{n-1} \leq \max(\tilde{x'})$.
Consequently, we have $\asctop(x) = \asctop(x') = \max(\tilde{x'}) = \max(\tilde{x})$.
If $x_{n-1}<x_n\leq \asctop(x')+1$, then $x_n \leq \max(\tilde{x'})+1$.
Recall that $\tilde{x} = \Add(\tilde{x'},x_n)$.
By the definition of the map $\Add$, we have $\max(\tilde{x}) = \max(\tilde{x'}) + 1 = \asctop(x') +1 = \asctop(x)$,  as desired.

Now we proceed to prove the second part of the lemma.
Let $y = \ha(x)$ and $z = \Add(w,x_n)$, where $w = \ha(x')$.
We need to show $y = z$.
From the definitions of $\ha$ and $\Add$,
we deduce that $y_2y_3\cdots y_{n+1} = z_2z_3\cdots z_{n+1} = \tilde{x}$.
Thus,  to prove $y=z$, we need only to show that $y_1 = z_1$.
Again from the definition of $\Add$, 
\begin{equation}\label{equ:z1}
	z_1=
	\begin{cases}
		w_1,     \quad &\text{if } 1\leq x_n\leq x_{n-1}, \\
		w_1 +1,  \quad&\text{if }  x_{n-1}< x_n \leq \asctop(x')+1.
	\end{cases}
\end{equation}
It is obvious that 
\begin{equation}\label{equ:asctop}
	\asctop(x)=
	\begin{cases}
		\asctop(x'),     \quad &\text{if } 1\leq x_n\leq x_{n-1}, \\
		\asctop(x') +1,  \quad&\text{if }  x_{n-1}< x_n \leq \asctop(x')+1.
	\end{cases}
\end{equation}
Note that $y_1 = \max(y) = \asctop(x)$.
Similarly, we have $w_1 = \max(w) = \asctop(x')$.
Hence, equation (\ref{equ:asctop}) can be rewritten as 
\begin{equation}\label{equ:y1}
	y_1=
	\begin{cases}
		w_1,     \quad &\text{if } 1\leq x_n\leq x_{n-1}, \\
		w_1 +1,  \quad&\text{if }  x_{n-1}< x_n \leq \asctop(x')+1.
	\end{cases}
\end{equation}
Comparing equations (\ref{equ:z1}) and 
(\ref{equ:y1}) gives  $y_1 = z_1$.
This completes the proof.
\qed

Lemma \ref{lem:Add} provides us an alternating definition of $\ha$, 
that is, $\ha(x) = \widetilde{1x}$, which can be proved by induction 
on the length of $x$.
The detailed proof is left to the reader.
We also need the following lemma to prove $\ha$ is a bijection between 
$\A_n$ and $\RA_{n+1}$.

% if and only if
%\begin{lemma}
%	For $n\geq 2$, let $x=x_1x_2\cdots x_n\in \RA_n$ be a revised ascent sequence and 
%	let $1\leq \ell\leq \max(x) +1$.
%	Then we have $\Add(x,\ell)\in \RA_{n+1}$.
%\end{lemma}

\begin{lemma}\label{lem:Add2}
	For $n\geq 2$, let $x=x_1x_2\cdots x_n$ be an endofunction and 
	let $1\leq \ell\leq \max(x) +1$.
%$x	\in \RA_n$ be a revised ascent sequence and 
%	let $1\leq \ell\leq \max(x) +1$.
	Then  $x	\in \RA_n$ if and only if $\Add(x,\ell)\in \RA_{n+1}$.
\end{lemma}

\pf
Let $y =y_1y_2\cdots y_{n+1}= \Add(x,\ell)$. 
If $x	\in \RA_n$, we proceed to show that $y\in \RA_{n+1}$, namely, 
 $y$ is a Cayley permutation such that $\Ascbot(y) = \nub(y)$.
We consider two cases. 
If $1\leq \ell \leq x_n$, then $y = \Add(x,\ell) = x\ell$.
Note that $x$ is a Cayley permutation with $\Ascbot(x) = \nub(x)$.
It follows that $\ell$ occurs in $x$, and thus $y$ is also a Cayley permutation.
Consequently, we deduce that  $\Ascbot(y) = \Ascbot(x) = \nub(x) = \nub(y)$, 
implying that $y\in \RA_{n+1}$.
If $x_n<\ell \leq \max(x)+1$, by the definition of $\Add$,
 $y$ contains each integer between 
$1$ and $\max(x)+1$.
Therefore, $y$ is a Cayley permutation.
Furthermore, for $1\leq i\leq n-1$, we have $i\in \nub(x)$ if and only if $i \in \nub(y)$.
Noting that the number $x_n$ occurs exactly once in $y$, 
it follows that $n\in \nub(y)$.
We claim that $n+1 \nin \nub(y)$. 
If not, the only possibility is that $\ell = x_n+1$ and the number $x_n$ occurs exactly once in $x$.
Then $x_n$ must serve as an ascent bottom in $x$,  which is impossible.
This proves the claim.
Thus, we have $\nub(y) = \nub(x)\cup \{n\}$.
By a similar argument as above, we have 
$\Ascbot(y) = \Ascbot(x)\cup \{n\}$.
Then $\Ascbot(y) = \nub(y)$ follows directly from the fact that $\Ascbot(x) = \nub(x)$.
Therefore, $y\in \RA_{n+1}$. 
%This completes the proof.

Conversely, if $y =\Add(x,\ell)\in \RA_{n+1}$, then
$x = \text{Rom}(y)$.
Combining the definition of the map Rom and the  fact that $y$ is
 a Cayley permutation, we have that $x$ is also a Cayley permutation.
Note that $\Ascbot(y) = \nub(y)$.
Using a similar argument as above,  we derive that if $y_{n+1}\leq y_n$,
then $\Ascbot(x) = \Ascbot(y) = \nub(y) = \nub(x)$;
and if $y_{n+1} > y_n$, then  $\Ascbot(x) = \Ascbot(y) \setminus \{n\} = \nub(y) \setminus \{n\} = \nub(x)$.
In both cases, we have  $\Ascbot(x)  = \nub(x)$,
thereby, $x\in \RA_n$.
\qed

\begin{theorem}\label{thm:hat}
	For $n\geq 1$, the map $\ha$ is a bijection between 
	$\A_n$ and $\RA_{n+1}$.
%	Furthermore, we have $\asctop(x) = \asctop(\ha(x))$ for any $x\in \A_n$.
  \end{theorem}

\pf
%We first show that $\ha$ is a map from $\A_n$ to $\MA_{n+1}$ by induction on $n$.
%Let $x = x
%The second part of the theorem is a direct consequence of the first part, along with Lemmas \ref{lem:Add} and \ref{lem:easy-prop}.
%Let us now consider the first part of the theorem.
Recall that $\ha$ is injective.
To prove  $\ha$ is a bijection between $\A_n$ and $\RA_{n+1}$, it suffices to show that the image 
$\ha(\A_n)$ is the set $\RA_{n+1}$.
We proceed to prove this result by induction on $n$ where $n=1$ is trivial. 
Assume the result for $n-1$. 
Then the result for $n$ can be verified in the following process.
%We need to prove this result for $n$.
\begin{align*}
	y\in \ha(\A_n)  &\iff y = \ha(x) \text{ for some } x\in \A_n \\
	&\iff y = \Add(\ha(x'),\ell) \text{ for some } x'\in \A_{n-1} \text{ and } 1\leq \ell \leq \asctop(x')+1 \\
	& \qquad\qquad\qquad\qquad\qquad\qquad (\mbox{by Lemma \ref{lem:Add}}) \\
	&\iff y =  \Add(y',\ell) \text{ for some } y' \in \RA_{n} \text{ and } 1\leq \ell \leq \max(y')+1 \\
	& \qquad\qquad\qquad\qquad\qquad\qquad (\mbox{by the induction hypothesis and Lemma \ref{lem:Add}})\\
	& \iff y \in \RA_{n+1}. \qquad\qquad\quad\;\; (\mbox{by Lemma \ref{lem:Add2}})
\end{align*}
Hence, $\ha(\A_n) = \RA_{n+1}$, 
as desired.
%This completes the proof.
\qed

From Theorem \ref{thm:hat}, we have $|\RA_{n}| = |\A_{n-1}|=F_{n-1}$.
Thus, revised ascent sequences constitute a Fishburn structure.
Let us end this section with a generating tree for revised ascent sequences.
%In the rest of this section, we will develop a generating tree for revised ascent sequences. 
The generating tree method was introduced by West \cite{West} and has since found widespread application in enumerative combinatorics.
Given $x \in \RA_{n}$ for $n\geq 2$, we say $y$ is a \emph{child} of $x$ if and only if $y =\Add(x,\ell)$ for $1\leq \ell \leq \max(x)+1$.
Based on Lemma \ref{lem:Add2}, this generating rule determines a generating tree for $\{\RA_n\}_{n\geq 2}$.
By the definition of $\Add$, the generating tree for $\{\RA_n\}_{n\geq 2}$
 is encoded by the following succession rules, where the two parameters $m$, $\ell$ track the maximum number and the last entry, respectively:
\begin{equation}\label{tree:RAS}
	\left\{
	\begin{aligned}
		\text{root} \colon \; & (1,1), \\[5pt]
		\text{rule} \colon \; & (m,\ell) \mapsto \{(m,i) \mid 1 \leq i \leq \ell\} \cup \{(m+1,i) \mid \ell < i \leq m+1\}.
	\end{aligned}
	\right.
\end{equation}

The root $(1,1)$ corresponds to the unique revised ascent sequence
of length $2$, which is the sequence $11$.
It is evident that ascent sequences  follow the same succession  rules  \cite{Cerb2406}, 
thereby reaffirming the existence of a bijection between $\A_n$ and $\RA_{n+1}$.

\section{Pattern avoidance in $\RA_n$}\label{sec:pattern}

In this section, we will investigate revised ascent sequences that avoid a single pattern.
Let us begin with some basic properties of revised ascent sequences.

%Moreover, one can verify that $$\ha(x) = \Add(1,x_1,x_2,\ldots,x_n).$$
%We leave the proof of this conclusion to the readers and  will not use it in this paper.

\begin{lemma}\label{lem:easy-prop}
	Let $n\geq 1$ and let $x=x_1x_2\cdots x_n$ be a revised ascent sequence. 
	Then we have 
	\begin{enumerate}
		\item [(1)] $x_1 = \max(x)$;
		\item [(2)] if $n\geq 2$, then $\max(x)$ occurs at least twice
		in $x$;
		\item [(3)] $\max(x) = \ascbot(x) = \asctop(x)$.
	\end{enumerate}
\end{lemma}

\pf
%Let $x$ be a revised ascent sequence in $\RA_n$ with $\max(x) = m$.
For convenience, we let $m = \max(x)$.
According to the definition of revised ascent sequences,  the leftmost copy of $m$ must serve as an ascent bottom in $x$, which 
can only happen if $x_1 = m$.
We now prove that $m$ occurs at least twice in $x$ if $n\geq 2$, which is trivial if $m=1$.
Assume that $m\geq 2$.
Since the leftmost copy of $m-1$ must serve as an ascent bottom in $x$,
the entry immediately following it must be $m$.
Consequently,  there are at least two copies of $m$ in $x$.

Note that $x$ is a Cayley permutation with $\Ascbot(x) = \nub(x)$.
Each integer between $1$ and $m$ serves as the leftmost entry of $x$ exactly once.
Consequently, we have $m = |\nub(x)| = \ascbot(x)$.
By the definitions of $\Ascbot(x)$ and $\Asctop(x)$, 
we see that the index $1$ belongs to both sets. 
Moreover, since $x_1= m$, it is apparent that for $k\geq 2$, 
$k\in \Ascbot(x)$ if and only if $k+1 \in \Asctop(x)$.
This implies that $\ascbot(x) = \asctop(x)$. 
\qed

%The following lemma will be used in our paper.

\begin{lemma}\label{lem:main}
	Let $\sigma = \sigma_1\sigma_2\cdots \sigma_k$ be a Cayley permutation. % with $\max(\sigma) = m$.
	Suppose that $\sigma_2 = m$ is the only occurrence of the maximum value in $\sigma$.
	Then we have $\RA_n(\sigma) = \RA_n(m\sigma)$.
\end{lemma}

\pf 
Let $x\in \RA_{n}$ be a revised ascent sequence.
If $x\in \RA_n(\sigma)$, then $x\in \RA_n(m\sigma)$.
Hence, $\RA_n(\sigma) \subseteq \RA_n(m\sigma)$.
To prove the opposite inclusion,  it is sufficient to prove that if $x$ contains $\sigma$, then $x$ contains $m\sigma$ too.
Suppose that $x_{i_1}x_{i_2}\cdots x_{i_k}$ is an occurrence of 
$\sigma$ in $x$.
Here we choose $x_{i_2}$ to be the leftmost largest such entry between 
$x_{i_1}$ and $x_{i_3}$.
Note that $\sigma_2$ is the only occurrence of the maximum value in $\sigma$.
We have either $x_{i_2}\geq x_{i_2+1}$ or $x_{i_2}$ is the rightmost entry in $x$.
In either case, we have $i_2 \nin \Ascbot(x)$.
Since $x$ is a revised ascent sequence, we have $\Ascbot(x) = \nub(x)$.
It follows that $i_2 \nin \nub(x)$. 
Due to the choice of $x_{i_2}$,  there is some $i_0<i_1$ such that 
$x_{i_0} = x_{i_2}$.
This implies that $x_{i_0}x_{i_1}x_{i_2}\cdots x_{i_k}$ is an occurrence of the pattern
$m\sigma$ in $x$, as desired.
\qed

\subsection{Patterns of length $2$ and pattern $212$}\label{subsection:12}

Let us begin with some simple patterns.
According to Lemma \ref{lem:easy-prop}, 
in each revised ascent sequence of length greater than  or equal to $2$, 
the maximum entry occurs at least twice.
Therefore, for $n\geq 2$,  there does not exist any revised ascent sequence of length $n$ that avoids the pattern $11$.
For $\sigma \in \{12,21\}$,  there is exactly one revised ascent sequence
of length $n$ that avoids the pattern $\sigma$, namely, the sequence $1^n$.
Here, $a^n$ denotes the sequence consisting of $n$ copies of $a$.
% ($n$ copies of $1$).
From Lemma \ref{lem:main}, we immediately obtain that $\RA_n(12) = \RA_n(212)$.
Therefore, the sequence $1^n$ is also the unique revised ascent sequence in $\RA_n(212)$.
%For the pattern $221$, we have the following proposition.

\subsection{Pattern $221$}\label{subsection:221}

\begin{proposition}
	For $n \geq 2$, we have $|\RA_n(221)| = n-1$.
\end{proposition}

\pf 
Let $x$ be a revised ascent sequence with $\max(x) = m$.
We claim that $x = m12\cdots (m-1)m^{n-m}$.
Indeed, since $x$ avoids the pattern $221$,  no entry smaller than $m$ can appear after the second copy of $m$ in $x$.
Furthermore, by Lemma \ref{lem:easy-prop}, we know that $x$ begins with $m$.
Thus,  $x$ must have the form $x = mym^k$, where each entry in $y$ is strictly smaller than $m$.
Note that $x$ is a revised ascent sequence.
We have
$\Ascbot(x) = \nub(x)$, meaning that the leftmost copy of each integer between 
$1$ and $m-1$ must serve as an ascent bottom.
This implies that there is a unique $m-1$ immediately preceding the second copy of $m$ in $x$.
Similarly, there is a unique $m-2$  immediately preceding the $m-1$, and so on.
Therefore, $y = 12\cdots m-1$ and $k = n-m$.
This proves our claim.
It is obvious that $m$ can range from $1$ to $n-1$, yielding that 
$|\RA_n(221)| = n-1$.
This completes the proof.
\qed

\subsection{Pattern $312$}\label{subsection:312}

\begin{proposition}
	For $n\geq 2$, we have $|\RA_n(312)| = 2^{n-2}$.
\end{proposition}

\pf
Let $x$ be a revised ascent sequence in $\RA_n(312)$ with $\max(x) = m$.
According to Lemma \ref{lem:easy-prop}, we have $x_1 = m$.
Since $x$ is a $312$-avoiding sequence, the entries strictly smaller than $m$ must appear in weakly 
decreasing order in $x$.
Furthermore, as $x$ is a Cayley permutation such that $\Ascbot(x) = \nub(x)$, it follows that  for each integer $k$ with $1\leq k \leq m-1$, the leftmost copy of $k$ in $x$ is followed immediately by at least one occurrence of $m$.
Hence, $x$ must be of the form 
\[x = m^{a_m}(m-1)m^{a_{m-1}}(m-1)^{b_{m-1}}(m-2)m^{a_{m-2}}(m-2)^{b_{m-2}}\cdots 1m^{a_1}1^{b_1},\]
where $a_i\geq 1$ for $1\leq i\leq m$ and $b_j\geq 0$ for $1\leq j\leq m-1$.
Moreover, we have 
\[\sum_{i=1}^{m}a_i  + \sum_{j=1}^{m-1}(b_j+1) = n.\]
One can verify that we have $\binom{n-1}{2m-2}$ ways to choose these
$a_i$ and $b_j$.
Therefore, we have \[|\RA_n(312)| = \sum_{m \geq 1}\binom{n-1}{2m-2} = 2^{n-2},\]
as desired.
\qed

\subsection{Pattern $122$}\label{subsection:122}

Given two positive integers $m$ and $n$ with $m\leq n$, 
let $[m,n] = \{m,m+1,\ldots, n\}$ and 
let $\overrightarrow{[m,n]}$ denote the sequence $m(m+1)\cdots n$.
To enumerate $\RA_{n}(122)$, we need the following description of the revised ascent sequences in $\RA_{n}(122)$.

%	Then $x$ must be the form $x = a_ia_i^{b_i}\overrightarrow{[a_{i-1},a_i]} a_{i-1}^{b_{i-1}}\overrightarrow{[a_{i-2},a_{i-1}]} \cdots a_2^{b_2}\overrightarrow{[a_1,a_2]}a_1^{b_1}\overrightarrow{[1,a_1]}$,

\begin{lemma}\label{lem:122}
	Given $n\geq 2$, let $x$ be a revised ascent sequence in $\RA_{n}(122)$ with $\max(x) = m$. 
	Then $x$ must be of the form $$x = a_1^{b_1+1}\overrightarrow{[a_{2},a_1]} a_{2}^{b_{2}}\overrightarrow{[a_{3},a_{2}]}a_3^{b_3} \cdots a_{i-1}^{b_{i-1}}\overrightarrow{[a_{i},a_{i-1}]}a_i^{b_i},$$
	where $m = a_1>a_2>\cdots > a_i = 1$
%	where $1=a_i<a_{i-1}<\cdots < a_1 = m$ 
	and $b_j\geq 0$ for $1\leq j \leq i$.
\end{lemma}

\pf
From Lemma \ref{lem:easy-prop}, $x$ begins with $m$.
Suppose that the leftmost entry which is not equal to $m$ is $a_2$.
Since $x \in \RA_{n}(122)$, each element of  $[a_2+1,m-1]$ occurs exactly once  in $x$.
Moreover, these elements must serve as ascent bottoms in $x$.
Hence, the entry $m-1$ must be immediately followed by a copy of $m$.
Note that there is at most one copy of $m$ to the right of the leftmost  $a_2$.
This implies that the if  $m-2 \in [a_2+1,m-1]$, then it is 
immediately followed by $m-1$.
By the same reasoning,  each entry $j \in [a_2+1,m-1]$ is immediately followed by $j+1$ in $x$.
Then the leftmost copy of $a_2$ must immediately followed by $a_2+1$.
To conclude, $x$ begins with $a_1^{b_1+1}\overrightarrow{[a_{2},a_1]}$ for some $b_1\geq 0$ and we assume 
$x = a_1^{b_1+1}\overrightarrow{[a_{2},a_1]}y$. 
If $a_2 \neq 1$, then $y$ is non-empty.
Suppose that the leftmost entry in $y$ which is not equal to $a_2$ is $a_3$.
Clearly, $a_3< a_2$. 
By a similar argument as above,  $y$ begins with $a_2^{b_2}\overrightarrow{[a_{3},a_2]}$ for some $b_2\geq 0$. 
This procedure can be repeated iteratively until we reach $a_{i-1}^{b_{i-1}}\overrightarrow{[a_{i},a_{i-1}]}$ with $a_i = 1$ for some $i$.
Finally, $x$ may end with some copies of $1$.
This completes the proof.
\qed

\begin{proposition}
	For $n\geq 2$,  we have $|\RA_n(122)| = 2^{n-2}$.
\end{proposition}

\pf
From Lemma \ref{lem:122}, each revised ascent sequence $x$ in $\RA_n(122)$ has the following form $$x = a_1^{b_1+1}\overrightarrow{[a_{2},a_1]} a_{2}^{b_{2}}\overrightarrow{[a_{3},a_{2}]} a_{3}^{b_{3}}\cdots  a_{i-1}^{b_{i-1}}\overrightarrow{[a_{i},a_{i-1}]}a_i^{b_i},$$
where $\max(x) = a_1>a_2>\cdots > a_i = 1$ and $b_j\geq 0$ for $1\leq j \leq i$.
Let $m = \max(x)$.
If $m = 1$ or $i = 1$, there is only one revised ascent sequence, namely, the sequence $1^n$.
For fixed $m\geq 2$ and $i\geq 2$, 
there are $\binom{m-2}{i-2}$ ways to choose the strictly decreasing sequence $m>a_2>a_3>\cdots > a_{i-1}>1$.
Since the length of $x$ is $n$,  we have 
%\[b_1+b_2+\cdots +b_i = n-m-i+1.\]
$b_1+b_2+\cdots +b_i = n-m-i+1$.
This implies that there are $\binom{n-m}{i-1}$ ways to choose these
$b_j$ for $1\leq j \leq i$.
Hence, we have 
\begin{align*}
	|\RA_n(122)| &= 1+\sum_{m\geq 2}\sum_{i\geq 2}\binom{m-2}{i-2}\binom{n-m}{i-1} \\
	&= 1+\sum_{m\geq 2}\binom{n-2}{m-1}\\
	&= 2^{n-2},
\end{align*}
as desired.
\qed

\subsection{Patterns $321$ and $231$} \label{subsection:321}

\begin{proposition}\label{prop:321}
	For $n \geq 1$, we have $|\RA_n(321)| = 2^{n-1}-n+1$.
\end{proposition}

\pf
Let $x$ be a revised ascent sequence in $\RA_n(321)$ with $\max(x) = m$.
If $m = 1$,  there is only one such revised ascent sequence, namely,  $x = 1^n$.
Now we consider $m\geq 2$.
According to Lemma \ref{lem:easy-prop}, we have $x_1 = m$.
Since $x$ avoids the pattern $321$, all entries strictly smaller than $m$ must appear
in weakly increasing order in $x$.
Furthermore,
each number $k \in [1,m-2]$ appears exactly once in $x$. 
Otherwise, the rightmost copy of $k$ would serve as an ascent bottom in $x$, 
contradicting to the fact $\Ascbot(x) = \nub(x)$.
Observe that, for each entry $k\in [1,m-2]$, it is either followed immediately by $k+1$ or by $m$.
Additionally, the leftmost $m-1$ must be followed immediately by 
$m$ in $x$. 
Finally, $x$ may end with some copies of $m-1$.
To conclude, $x$ must be of the following form:
\[x = m^{a_0+1}1m^{a_1}2m^{a_2}\cdots (m-2)m^{a_{m-2}}(m-1)m^{a_{m-1}+1}(m-1)^{a_m},\]
where $a_i\geq 0$ for $0\leq i\leq m$.
Moreover, we have 
\[\sum_{i=0}^{m}a_i = n-m-1\] as the length of $x$ is $n$.
Since there are  $\binom{n-1}{m}$ possible choices for these $a_i$, 
we arrive at the result
\[|\RA_n(321)| = 1+\sum_{m\geq 2}\binom{n-1}{m} = 2^{n-1}-n+1.\]
This completes the proof.
\qed

\begin{proposition}\label{prop:231}
	For $n \geq 1$, we have $|\RA_n(231)| = 2^{n-1}-n+1$.
\end{proposition}

\pf
We proceed to prove this conclusion by showing that
$\RA_n(231) = \RA_n(321)$.
Assume that $x\in \RA_n(321)$.
Recall that all entries strictly smaller than $\max(x)$ appear
in weakly increasing order in $x$.
It follows that $x$ avoids the pattern $231$.
Thus, we have $\RA_n(321) \subseteq \RA_n(231)$.

Conversely, assume that $x\in \RA_n(231)$ with $\max(x) = m$.
To prove that $x\in \RA_n(321)$, it suffices to show
that all entries in the set $[1,m-1]$ are weakly increasing in $x$.
Suppose that there exist indices $i<j$ such that  $x_j<x_i\leq m-1$.
Choose $x_i$ to be the leftmost such entry.
Since $x$ is a revised ascent sequence, 
we have $x_i<x_{i+1}$.
This implies that the subsequence $x_ix_{i+1}x_j$ forms an occurrence of the pattern $231$ in $x$, which contradicts the assumption that $x\in \RA_n(231)$.
Therefore, we conclude that all entries
in the set $[1,m-1]$ must appear in weakly increasing order in $x$, proving that 
 $x\in \RA_n(321)$.
 Hence, we have  $\RA_n(231) \subseteq \RA_n(321)$.
 
Then we conclude that $\RA_n(231) = \RA_n(321)$.
Combining this with Proposition \ref{prop:321} gives 
$|\RA_n(231)| = |\RA_n(321)| =2^{n-1}-n+1$.
\qed

The following result follows directly from Lemma \ref{lem:main}
and Proposition \ref{prop:231}.
\begin{proposition}
	For $n \geq 1$, we have $|\RA_n(3231)| = 2^{n-1}-n+1$.
\end{proposition}

\subsection{Pattern $213$}\label{subsection:213}

Let $x$ be a sequence of positive integers.
The \emph{standardization} of $x$, denoted by $\std(x)$, is defined to be the sequence which is obtained from $x$ by 
relabeling the sequence according to their relative order.
More specifically, the smallest integers in $x$ are replaced with $1$, the next smallest with $2$, and so on.
As an example, we have $\std(7424326) = 5313214$.
Given $n \geq 1$, let $f_n = |\RA_n(213)|$.
Then we have the following recurrence for $f_n$.

\begin{lemma}\label{lem:213}
	For $n\geq 1$, we have 
	\begin{equation}\label{equ:213}
		f_n = 1+\sum_{k= 1}^{n-2}f_k\sum_{j=2}^{n-k}f_j.
	\end{equation}
%\[	f_n = 1+\sum_{k= 1}^{n-2}f_k\sum_{j=2}^{n-k}f_j.\]
%with the initial value $f_1=1$.
\end{lemma}

 \pf 
 There is exactly one revised ascent sequence of length $1$ and one of length $2$ that avoids the pattern $213$,
 namely, $1$ and $11$, respectively.
 Hence, equation (\ref{equ:213}) holds for $n=1,2$.
 Now we assume that $n\geq 3$.
 Let $x \in \RA_{n}(213)$ with $\max(x) = m$.
There is exactly one revised ascent sequence in $\RA_{n}(213)$ satisfying $\max(x) = 1$,
 namely, the sequence $1^n$.
 And now we assume that $m\geq 2$. 
 We decompose the sequence $x$ into the form $x = L1R$ according to 
the leftmost $1$ in $x$.
Let $k$ be the length of $L$.
Since $x$ avoids the pattern $213$, 
all entries in $L$ are greater than or equal to the entries in $R$.

%Then for any entry $a$ in $L$ and $b$ in $R$, we have 
%$a\geq b$.
%Otherwise, if there is some $a$ in $L$ and some $b$ in $L$ such that $a<b$.
%By our decomposition, it is obvious that each entry in $L$ is strictly greater than $1$.
%This implies that $a1b$ forms an occurrence of $213$ in $x$, 
%a contradiction to the fact that $x$ is $213$-avoiding.
%Hence, all entries in $L$ are greater than or equal to the entries in $R$.

Let $s =\min(L)$.
%Then all entries in $R$ are smaller than or equal to $s$.
Since the leftmost $s-1$ must serve as an ascent bottom in $R$, 
we deduce that the entry $s$ appears in $R$.
Notice that no other $1$'s except for the leftmost $1$ can serve as an ascent bottom
in $x$.
Thus, if there are any other $1$'s in $x$, they can only appear at the end of $R$.
%Let $i$ denote the number of consecutive $1$s that appear at the end of $R$.
Then $x$ can be further written as $x = L1R'1^i$ for some $i\geq 0$.
Note that $R'$ contains at least one copy of $s$.
We have $0\leq i\leq n-k-2$.
It is easily seen that 
$\std(L) \in \RA_{k}(213)$ and $\std(sR')\in \RA_{n-k-i}(213)$.
Conversely, one can recover a unique $x \in \RA_n(213)$ from 
a revised ascent sequence in $\RA_{k}(213)$ and a revised ascent sequence in 
$\RA_{n-k-i}(213)$.
Consequently, we arrive at the desired equation (\ref{equ:213}), with the summation variable
$j$ indicating the length of $sR'$.
\qed

\begin{proposition}
	Let $C(t) = \sum_{n\geq 0}C_nt^n = \frac{1-\sqrt{1-4t}}{2t}$ be the generating function of the Catalan number $C_n$. 
	Then we have 
	\begin{equation}\label{ogf:213}
		\RA_{213}(t) = t\,C(t) = \frac{1-\sqrt{1-4t}}{2}.
	\end{equation}
	Therefore, $f_n = |\RA_{n}(213)| = C_{n-1}$.
\end{proposition}

\pf
Turning the recurrence relation (\ref{equ:213}) into generating function yields
\[\RA_{213}(t) = \frac{t}{1-t} + \RA_{213}(t) \frac{\RA_{213}(t) - t}{1-t}.\]
Solving the above equation, we obtain (\ref{ogf:213}), as desired.
\qed

\subsection{Patterns $211$ and $121$}\label{subsection:211}

%Let $\MA_n(122)$ denote the set of modified ascent sequences in $\MA_n$ that avoid
%the pattern $122$.
%\begin{lemma}[Cerbai, \cite{Cerbai23a}]
%	For $n\geq 1$, we have 
%	\[|\MA_n(122)| = \sum_{k=1}^{n}k^{n-k}.\]
%\end{lemma}

\begin{proposition}\label{prop:211}
		For $n\geq 2$, we have 
		\[|\RA_n(211)| = \sum_{k=1}^{n-1}k^{n-k-1}.\]
%		\begin{equation}\label{equ:211}
%			|\RA_n(211)| = \sum_{k=1}^{n-1}k^{n-k-1}.
%		\end{equation}
\end{proposition}

 \pf 
 Let $x$ be a revised ascent sequence in $\RA_n(211)$ with $\max(x) = m$.
According to Lemma \ref{lem:easy-prop}, we have that $x_1 = m$.
 Since $x$ avoids the pattern $211$,  each integer between $1$ and $m-1$ appears exactly once in $x$.
 Note that $\nub(x) = \Ascbot(x)$, 
 implying that these integers must serve as ascent bottoms 
 in $x$.
 Suppose that $x$ contains $k+1$ ($1\leq k\leq n-1$) copies of $m$.
 Then $x$ can be decomposed into the following form
\[x = mB_1mB_2\cdots mB_km,\]
where the entries in each block $B_i$ are strictly less than $m$.
Moreover, each $B_i$ is a strictly increasing sequence (possibly empty).
Now each integer in the set $[1,m-1]$ can be assigned to one of the $k$ blocks.
Since the length of $x$ is $n$, we deduce $m+k = n$.
Thus, the number of revised ascent sequences in $\RA_n(211)$ with $k+1$ copies of 
maximum number is $k^{m-1}$, namely, $k^{n-k-1}$.
Summing over $k$, we get 
\[|\RA_n(211)| = \sum_{k=1}^{n-1}k^{n-k-1},\]
as desired.
\qed

Let $\MA_n(122)$ denote the set of modified ascent sequences in $\MA_n$ that avoid the pattern $122$.
By comparing Proposition \ref{prop:211} with the enumeration
formula for $\MA_n(122)$ due to Cerbai \cite{Cera}, 
we arrive at
 $|\RA_{n+1}(211)| = |\MA_{n}(122)|$.
 We construct a bijection $\phi$ between $\RA_{n+1}(211)$ and $\MA_{n}(122)$ as follows.
 Given $x\in \RA_{n+1}(211)$,  
 \begin{enumerate}
 	\item increment each entry of  $x$ by 1;
 	\item replace each copy of the largest entry in the modified sequence with a $1$ to ensure that the sequence qualifies as a Cayley permutation;
 	\item remove the last entry of the sequence.
 \end{enumerate}
The sequence resulting from these transformations is defined as $\phi(x)$.
% we first increment each entry of  $x$ by 1. 
% Next, we replace each copy of the largest entries in the modified sequence with a $1$ to ensure that the sequence qualifies as a Cayley permutation. 
% Finally, we remove the last entry of the sequence. 
%  The resulting sequence obtained through these transformations is defined as $\phi(x)$.
For example, given $x = 6463612656 \in \RA_{10}(211)$, then we have 
$\phi(x) = 151412316$ with the following detailed process:
\begin{align*}
	x = 6463612656 \rightarrow 7574723767 \rightarrow 1514123161 \rightarrow 151412316 = \phi(x).
\end{align*}
Based on the characterization of the set $\MA_{n}(122)$ due to 
Cerbai \cite{Cera}, it is easy to verity that $\phi$ is a bijection 
between $\RA_{n+1}(211)$ and $\MA_{n}(122)$.
 Let $\mathcal{R}_n$ denote the set of partitions
of $[n]$ whose minima of all blocks form an interval.
Cerbai \cite{Cera} also constructed a bijection between $\MA_{n}(122)$ and 
$\mathcal{R}_n$.
It yields that $|\RA_{n+1}(211)| = |\mathcal{R}_n|$.
By combining $\phi$ and the bijection between $\MA_{n}(122)$ and 
$\mathcal{R}_n$ constructed by Cerbai, 
one can also obtain a bijection between  $\RA_{n+1}(211)$ and 
$\mathcal{R}_n$.

 %Now we consider the pattern $121$.
 \begin{proposition}\label{prop:121}
 			For $n\geq 2$, we have 
 	\[|\RA_n(121)| = \sum_{k=1}^{n-1}k^{n-k-1}.\]
 \end{proposition}

\pf
We proceed to prove this conclusion by demonstrating that $\RA_n(121) = \RA_n(211)$.
Let $x$ be a revised ascent sequence in $\RA_n$, namely, 
$x$ is a Cayley permutation with $\Ascbot(x) = \nub(x)$.
Suppose that $x$ contains the pattern $121$.
By Lemma \ref{lem:main},
$x$ also contains the pattern $2121$.
Consequently, $x$ contains the pattern $211$.

Conversely, suppose $x$ contains the pattern $211$ and let $x_ix_jx_k$ be
an occurrence of $211$ in $x$.
Let $s$ be the index such that $s\in \nub(x)$ and $x_s = x_j$.
From Lemma \ref{lem:main}, we have  $x_1 = \max(x)$.
It follows that $1<s\leq j$.
Note that $\nub(x) = \Ascbot(x)$.
We have $s\in \Ascbot(x)$, which implies that $x_s<x_{s+1}$.
Hence, $x_sx_{s+1}x_k$ forms an occurrence of the pattern  $121$ in $x$.
To conclude, we have shown that $x$ contains the pattern $121$ if and only if $x$ contains the pattern 
$211$, namely, $\RA_n(121) = \RA_n(211)$.
Then the result follows directly from Proposition \ref{prop:211}, 
completing the proof.
\qed

Combining Lemma \ref{lem:main}
and Proposition \ref{prop:121}, we directly obtain the following proposition.
%The following proposition follows directly from Lemma \ref{lem:main}
%and Proposition \ref{prop:121}.
\begin{proposition}
		For $n\geq 2$, we have 
	\[|\RA_n(2121)| = \sum_{k=1}^{n-1}k^{n-k-1}.\]
\end{proposition}

%Fix $i$ and $k$, we choose the leftmost such $x_j$.
%To show that $x$ contains the pattern $121$, 
%we consider two cases.
%If $j \in \nub(x) = \ascbot(x)$, then $x_j<x_{j+1}$.
%This implies that $x_jx_{j+1}x_k$ forms an occurrence of $121$ in $x$.
%If $j\nin \nub(x)$, by our choice of $x_j$, there is some index $s<i$ such that $s\in \nub(x)$ and $x_s = x_j$.
%Then $x_sx_{s+1}x_

\subsection{Pattern $112$}\label{subsection:112}

Given two positive integers $n$ and $m$, let $\RA_{n,m}(112)$ denote the 
set of revised ascent sequences in $\RA_n(112)$ which have the maximum value $m$.
Let $S(n,m)$ denote the $(n,m)$-th Stirling number of the second kind.
It is well known that $S(n,m)$ is equal to the number of partitions of $[n]$ into
$m$ blocks.
We have the following proposition for $\RA_{n,m}(112)$.

\begin{proposition}
	For $n,m\geq 1$, we have $|\RA_{n,m}(112)| = S(n-1,n-m)$.
	Consequently, $|\RA_{n}(112)|$ %= b(n-1)$, where $b(n-1)$ 
	is equal to the $(n-1)$-th Bell number.
\end{proposition}

\pf
The conclusion is trivial when $n=1$ or $m=1$.
Now we assume that  $n,m\geq 2$.
Given $x \in \RA_{n,m}(112)$, let $a_i$ ($1\leq i\leq m$) denote  the number of occurrences of $i$ in $x$.
Then we have $\sum_{i=1}^{m}a_i = n$.
Based on the fact that $x$ is a Cayley permutation, we have $a_i\geq 1$ for $1\leq i\leq m$.
Moreover, by Lemma \ref{lem:easy-prop}, we have $a_m \geq 2$.

%To construct such a revised ascent sequence $x \in \RA_{n,m}(112)$, 
%we begin by positioning $a_m$ occurrences of the value $m$.
%Following this, we insert $a_{m-1}$ occurrences of $m-1$.
%The leftmost instance of $m-1$ must serve as an ascent bottom, and therefore, there are 
%$a_m-1$ ways to place it.
%Additionally, to ensure that $x$ avoids the pattern $112$, any remaining 
%$(m-1)$'s (if exist) must be placed at the end of the sequence.
%Next, we insert $a_{m-2}$ occurrences of $m-2$.
%To avoid the pattern $112$, all but the leftmost $m-2$ must  also be positioned at the end of the sequence.
%When placing the leftmost $m-2$, it is necessary for it to serve as an ascent bottom, yet it cannot directly follow the leftmost $m-1$.
%Consequently, there are $a_m+a_{m-1}-2$ potential locations for this element.
%Continuing in this manner, %assume that we have placed $m,m-1,\ldots,i+1$.
%when placing $i$, the leftmost $i$ has $a_m+a_{m-1}+\cdots+a_{i+1}-(m-i)$ possible positions, while any remaining $i$ can only be placed at the end of the sequence.
%Hence, we have 

To construct such a revised ascent sequence $x \in \RA_{n,m}(112)$, we proceed as follows:
\begin{itemize}
	\item We begin by positioning $a_m$ occurrences of the value $m$.
	\item Next, we insert $a_{m-1}$ occurrences of $m-1$.
	The leftmost $m-1$ must serve as an ascent bottom, and therefore, there are $a_m-1$ ways to place it.
	Additionally, to ensure that $x$ avoids the pattern $112$, any remaining 
	$(m-1)$'s (if exist) must be placed at the end of the sequence.
	\item Then we insert $a_{m-2}$ occurrences of $m-2$.
	Similarly, 
	to avoid the pattern $112$,  all but the leftmost $m-2$ must  also be positioned at the end of the sequence.
	When placing the leftmost $m-2$, it is necessary for it to serve as an ascent bottom, yet it cannot directly follow the leftmost $m-1$.
	Consequently, there are $a_m+a_{m-1}-2$ potential locations for this element.
	\item Continuing in this manner, %assume that we have placed $m,m-1,\ldots,i+1$.
	when inserting $i$, the leftmost $i$ has $a_m+a_{m-1}+\cdots+a_{i+1}-(m-i)$ possible positions, while any remaining $i$'s can only be placed at the end of the sequence.
\end{itemize}
Hence, we have 
\begin{align*}
	|\RA_{n,m}(112)| &= \sum_{a_1+a_2+\cdots +a_m = n}\prod_{i=1}^{m-1}(a_m+a_{m-1}+\cdots+a_{i+1}-m+i) \\
	&= \sum_{1\leq b_{m-1}\leq \cdots \leq b_2\leq b_{1}\leq n-m}b_1b_2\cdots b_{m-1}  \quad (\text{by setting } b_i = a_m+a_{m-1}+\cdots+a_{i+1}-m+i)\\
	& = \sum_{\substack{c_1+c_2+\cdots + c_{n-m} = m-1\\ c_i\geq 0}}1^{c_1}2^{c_2}\cdots (n-m)^{c_{n-m}}
	\quad (\text{by setting $c_i$ to be the occurrences} \\
	&\qquad\qquad\qquad\qquad\qquad\qquad\qquad	\qquad\qquad\text{of $i$ in the sequence $b_1b_2\cdots b_{m-1}$}) \\
	& = S(n-1,n-m),
\end{align*}
where the last step follows by comparing the coefficients of $x^n$ on both sides of the following 
generating function \cite{Stanley}
\[\sum_{n}S(n,k)x^n = \frac{x^k}{(1-x)(1-2x)\cdots (1-kx)}. \]
This completes the proof.
\qed

In the proof above, if we let $d_i = c_i+1$ for $1\leq i\leq n-m$, then 
\[|\RA_{n,m}(112)| = \sum_{\substack{d_1+d_2+\cdots + d_{n-m} = n-1\\ d_i\geq 1}}1^{d_1-1}2^{d_2-1}\cdots (n-m)^{d_{n-m}-1}.\]
This means that each revised ascent sequence $x\in \RA_{n,m}(112)$ can be associated with 
a composition $d_1+d_2+\cdots +d_{n-m} = n-1$ and  exactly $1^{d_1-1}2^{d_2-1}\cdots (n-m)^{d_{n-m}-1}$ revised ascent sequences $x \in \RA_{n,m}(112)$  are associated with this composition.
On the other hand, we can also associate with each partition $\pi$ of $[n-1]$ into $n-m$ blocks a composition $d_1+d_2+\cdots +d_{n-m} = n-1$ by defining 
$d_{i+1}+d_{i+2}+\cdots + d_{n-m}$ to be the least $r$ such that when $1,2,\ldots,r$ are removed from $\pi$, the resulting partition has $i$ blocks.
And there are exactly $1^{d_1-1}2^{d_2-1}\cdots (n-m)^{d_{n-m}-1}$ partitions $\pi$ 
which are associated with this composition.
See Exercise $44$ of Chapter $1$ in \cite{Stanley} for details.
Hence, one can provide a combinatorial proof for $|\RA_{n,m}(112)| = S(n-1,n-m)$ by 
constructing a bijection between $\RA_{n,m}(112)$ and partitions of $[n-1]$ into $n-m$ blocks.
We omit the detailed specifics here, leaving them to the interested reader.

\subsection{Pattern $123$}\label{subsection:123}

In this section, we shall enumerate $\RA_n(123)$ by constructing its generating tree.
Recall that we have constructed a generating tree for revised ascent sequences $\{\RA_n\}_{n\geq 2}$ in the end of Section \ref{sec:bijection} with the succession rules as shown in (\ref{tree:RAS}).
Our generating tree for $\{\RA_n(123)\}_{n\geq 2}$, denoted by 
$\T(123)$,
 can be obtained from the one 
for $\{\RA_n\}_{n\geq 2}$ by restricting the vertices on the $123$-avoiding revised ascent sequences.
That is, given a revised ascent sequence $x \in \RA_n(123)$ with $n\geq 2$,
the set of children of $x$ is given by 
\[\{y \mid y = \Add(x,\ell) \text{ for $1\leq \ell \leq \max(x)+1$  and } y \in \RA_{n+1}(123)\}.\] 
To obtain the  succession rules for  $\T(123)$, we need a new labeling 
scheme.
For $x\in \RA_n(123)$,  define 
\[g(x) = \begin{cases} 
	1, & \text{if $x$ avoids $12$;} \\
	\min\{x_j \mid \text{there exists  $i$ such that $i<j$ and $x_i<x_j$}\}, & \text{otherwise.}
\end{cases} 
\]
For example, given $x = 5453423133\in \RA_{10}(123)$, we have 
$g(x) = 3$.
Indeed, for any revised ascent sequence $x\in \RA_n(123)$, $g(x) = 1$ if $x$ avoids the pattern $12$;
otherwise, $g(x)$ is the smallest entry among the larger elements of $12$ patterns of $x$.
Now we associate a label $(g(x),\ell(x))$ to $x$, where 
$\ell(x) = x_n$.
As an example, the label of the  revised ascent sequence $5453423133$ is $(3,3)$.

\begin{lemma}\label{lem:123}
The generating tree $\T(123)$  is encoded by the following succession rules:
\begin{equation}\label{tree:123}
	\left\{
	\begin{aligned}
		\mathrm{root} \colon\; & (1,1), \\[5pt]
		\mathrm{rule} \colon \;& (g,\ell) \mapsto \{(g,1)\} \cup \{(i,i) \mid 2 \leq i \leq g + \delta_{1\ell}\},
	\end{aligned}
	\right.
\end{equation}
where $\delta_{1\ell}$ is the Kronecker delta function, which equals $1$ if $\ell = 1$ and $0$ otherwise.
\end{lemma}

\pf 
It is apparent that 
the root $(1,1)$ corresponds to the unique revised ascent sequence
of length $2$, which is the sequence $11$.
Now we focus on the rule for the generating tree.
For $n\geq 2$, let $x$ be a revised ascent sequence in $\RA_n(123)$ with the label $(g,\ell)$.
That is, $g(x) = g$ and $x_n = \ell$.
Given $1\leq i\leq \max(x)+1$, let $y = \Add(x,i)$.
By applying Lemma \ref{lem:Add2}, we obtain $y\in \RA_{n+1}$.
Hence, to determine whether $y$ is a child of $x$, it suffices to 
determine whether $y$ avoids the pattern $123$.
Now we consider two cases.

\noindent {\bf Case (\rmnum{1}).} $\ell = 1$.\\
If $x$ avoids the pattern $12$, then $x = 1^n$ and $g = 1$.
In this case, 
 $y$ avoids the pattern $123$ if and only if $1\leq i\leq 2$.
 That is, the children of $x$ are $\{1^{n+1}, 2^{n-1}12\}$,  whose labels are the set
 $\{(1,1),(2,2)\}$.
One can easily verify that it follows the rule stated in (\ref{tree:123}).
Now we assume that $x$ contains the pattern $12$.
By the definitions of $\Add$ and $g(x)$,  $y$ avoids the pattern  $123$ if and only if
$1\leq i\leq g+1$.
If $i = 1$, the label of $y$ is $(g,1)$.
If $2\leq i\leq g+1$, then $g(y) = i$, 
implying that the label of $y$ is $(i,i)$.
Thus, the set of labels of the children of $x$ is $\{(g,1)\} \cup \{(i,i) \mid 2 \leq i \leq g + 1\}$, 
which also satisfies the rule stated in (\ref{tree:123}).

\noindent {\bf Case (\rmnum{2}).} $\ell \geq 2$.\\
Since $x$ avoids the pattern  $123$,  we deduce that  $g = \ell$.
By the definitions of $\Add$ and $g(x)$,  it is easily seen that if $1\leq i\leq g$, 
then $y$ avoids the pattern  $123$;
if $i>g=\ell$, then 
$1\ell i$ would form a subsequence of $y$, which is an occurrence of the pattern $123$.
Hence, $y$ avoids the pattern  $123$ if and only if $1\leq i\leq g$.
By a similar argument as Case (\rmnum{1}), we derive that the set of labels of the children of $x$ is $\{(g,1)\} \cup \{(i,i) \mid 2 \leq i \leq g\}$, 
as desired.
\qed

Given $n\geq 2$ and $x \in \RA_n(123)$, let 
$\T_x(123)$ be the subtree of $\T(123)$ consisting of 
the revised ascent sequence $x$ as the root and 
its descendants in $\T(123)$.
We define $F_x(t)$ to be the generating function for 
the number of vertices at level $n$ in the subtree $\T_x(123)$, 
where the root $x$ of the subtree $\T_x(123)$ is at level $1$.
That is, 
\[F_x(t) = \sum_{n\geq 1}|\{\text{vertices at level $n$ in the subtree $\T_x(123)$}\}|  t^n.\]
Notice that $F_x(t)$ is uniquely determined by the label $(g,\ell)$ of $x$. 
This means that  $F_x(t) = F_y(t)$ if and only if the two revised ascent sequences $x$ and $y$ have the same labels. 
Hence, we often write $F_{(g,\ell)}(t)$ for $F_x(t)$.
% if the label of $x$ is $(g,\ell)$.

For convenience, we let $A_g(t) = F_{(g,1)}(t)$ for $g\geq 1$ and 
$B_{\ell}(t) = F_{(\ell,\ell)}(t)$ for $\ell\geq 2$.
From the proof of Lemma \ref{lem:123}, 
we see that $g = \ell$ if $\ell\geq 2$.
Then the rule (\ref{tree:123}) can be translated into the following system of equations by distinguishing whether $\ell = 1$.
\begin{align}
	A_g(t) &= t+tA_g(t) + t\sum_{i=2}^{g+1}B_i(t) \quad (g\geq 1),\label{equ:Ag}\\ 
	B_{\ell}(t) &= t+tA_{\ell}(t) + t\sum_{i=2}^{\ell}B_i(t) \quad(\ell \geq 2).  \label{equ:Bg}
\end{align}
%Let $g = \ell$ in (\ref{equ:Ag}). 
%Subtracting equation (\ref{equ:Bg}) from (\ref{equ:Ag}) gives the following equation:
%\begin{equation}
%	A_{\ell}(t) = B_{\ell}(t) + tB_{\ell+1}(t)  \quad(\ell \geq 2). 
%\end{equation}

In order to solve the above system of equations, we define 
\begin{align*}
	A(t,u) &= \sum_{g\geq 1}A_g(t)u^{g-1},  \text{ and }\\
		B(t,u)& = \sum_{\ell\geq 2}B_{\ell}(t)u^{\ell-2}.
\end{align*}

\begin{proposition}
	We have 
	\begin{equation}\label{equ:123}
		\RA_{123}(t) = \sum_{n\geq 1}|\RA_n(123)|t^n = \frac{t^{2} -1+ \sqrt{t^{4} + 2 \, t^{2} - 4 \, t + 1} }{2 \, {\left(t - 1\right)}}. 
	\end{equation}
\end{proposition}

\pf
Multiplying both sides of equation (\ref{equ:Ag})  by $u^{g-1}$ and summing over all  $g\geq 1$, then simplifying yields
\begin{equation}\label{equ:A}
	A(t,u) = \frac{t}{1-u} + tA(t,u)+\frac{tB(t,u)}{1-u}.
\end{equation}
Similarly, multiplying both sides of equation (\ref{equ:Bg})  by $u^{\ell-1}$ and summing over all  $\ell\geq 2$, then simplifying yields
\begin{equation}\label{equ:B}
	uB(t,u) = \frac{ut}{1-u} + t\big(A(t,u)- A(t,0)\big)+\frac{utB(t,u)}{1-u}.
\end{equation}
Combining equations (\ref{equ:A}) and (\ref{equ:B}) and eliminating $B(t,u)$, we derive the following functional equation for $A(t,u)$.
\begin{equation}\label{equ:A0}
	\big( (t-1)u^2 + (t^2-2t+1)u -t^2  \big)A(t,u) =ut - t^2A(t,0).
\end{equation}

We proceed to solve equation (\ref{equ:A0}) by using the kernel method. 
Let $u_0 = u_0(t)$ satisfy 
\[(t-1)u^2 + (t^2-2t+1)u -t^2 = 0,\]
namely, 
\[u_0 = \frac{-t^{2} + 2 \, t - 1+\sqrt{t^{4} + 2 \, t^{2} - 4 \, t + 1} }{2 \, {\left(t - 1\right)}}.\]
(Note that there are two solutions for $u_0$, we choose the one by the condition $A(0,0) = 0$.)
Letting $u= u_0$ in (\ref{equ:A0}) and solving for $A(t,0)$ gives 
\begin{equation}
	A(t,0) = \frac{-t^{2} + 2 \, t -1+ \sqrt{t^{4} + 2 \, t^{2} - 4 \, t + 1} }{2 \, t{\left(t - 1\right)}}.
\end{equation}
Then (\ref{equ:123}) follows directly from the fact 
$\RA_{123}(t) = t+tA(t,0)$.
\qed 

From (\ref{equ:123}), we derive the sequence 
\[
(|\RA_{n}(123)|)_{n\geq 2} = (1,2,4,9,22,57,154,429,1223,\cdots),
\]
which appears as A105633 in the OEIS \cite{oeis}.
This sequence also counts 
the Dyck paths of semilength $n$ with no {\bf uudu} (or with no {\bf uduu}).

\subsection{Pattern $132$}\label{subsection:132}

For $n\geq 1$, let $g_n = |\RA_n(132)|$. 
We proceed to enumerate $g_n$ by constructing a system of equations for $g_n$, $r_n$ and $s_n$, where 
\begin{align*}
	r_n &= |\{x\in\RA_n(132)\mid x_n = \max(x) \}| \text{ and }\\
	s_n &= |\{x\in\RA_n(132)\mid x_{n-1}<\max(x) \text{ and }x_n = \max(x) \}|.
\end{align*}
%$g_n = |\RA_n(132)|$ and let $r_n = |\{x\in\RA_n(132)\mid x_n = \max(x) \}|$.
%For $n\geq 2$, let $s_n = |\{x\in\RA_n(132)\mid x_{n-1}<\max(x) \text{ and }x_n = \max(x) \}|$
For convenience, we set $g_0 = r_0 = s_0 = 0$  and $s_1 = 0$.
%We proceed to construct a system of equations for $g_n$, $r_n$ and $s_n$.

\begin{lemma}
	For $n\geq 2$, we have 
	\begin{equation}\label{equ:rs}
		r_n = r_{n-1}+s_n.
	\end{equation}
\end{lemma}

\pf 
Given a revised ascent sequence $x\in \RA_n(132)$ with $x_n = \max(x)$,
we have either $x_{n-1}<\max(x)$ or $x_{n-1} = \max(x)$.
In the first case, it is counted by $s_n$.
For the latter case, by removing the last entry of $x$, we obtain
 a revised ascent sequence in $\RA_{n-1}(132)$
 that ends with its maximum entry.
 Hence, it is counted by $r_{n-1}$.
Then equation (\ref{equ:rs}) follows.
 \qed

\begin{lemma}
	For $n\geq 1$, 
	\begin{equation}\label{equ:gr}
		g_n = 1+\sum_{i = 3}^{n}(r_i -1)g_{n+1-i}.
	\end{equation}
\end{lemma}

\pf
Let $x$ be a revised ascent sequence in $\RA_n(132)$ with $\max(x) = m$.
If $m = 1$, we have $x = 1^n$.
Now we assume that $m \geq 2$.
The sequence $x$ can be  decomposed uniquely as $x = LmR$ according to the rightmost $m$ of $x$.
Since $x$ avoids the pattern $132$, we derive that each entry in $L$ is
not smaller than any entry in $R$.

Let $m' = \min(L)$ and let $i$ be the length of $Lm$.
We  define $L'$ and $R'$ as the  standardizations of 
$Lm$ and $m'R$, respectively.
That is, $L' = \std(Lm)$ and $R' = \std(m'R)$.
Based on the fact that $\Ascbot(x) = \nub(x)$, 
one can easily verify that $L'$ is a revised ascent sequence in $\RA_i(132)$ ending with its maximum entry, while $R'$ is a revised ascent sequence in $\RA_{n+1-i}(132)$.
Furthermore, since the leftmost $m-1$  serves as an ascent bottom in $x$, 
it must appear in $L$.
Note that $x_1 =m$.
It follows that $i\geq 3$ and $L'$ cannot be composed entirely of $1$'s.
Thus, we have $r_i - 1$ ways to choose $L'$ and $g_{n+1-i}$ ways to choose 
$R'$.
Conversely, one can easily recover the unique revised ascent sequence $x\in \RA_n(132)$ from $L'$ and $R'$.
Hence, equation (\ref{equ:gr}) follows.
\qed

\begin{lemma}
	For $n\geq 2$, we have 
	\begin{equation}\label{equ:gs}
		s_n = \sum_{k\geq 1}g_{n-1-k}.
	\end{equation}
\end{lemma}

\pf 
Let $x$ be a revised ascent sequence in $\RA_n(132)$ such that $x_{n-1}<x_n = \max(x)$. 
For convenience, let $m = \max(x)$ and $k = x_{n-1}$.
Then we have  $k<m$.
We claim that $x$ is of the form $x = y12\cdots k m$.
If $k=1$, then there is nothing to prove.
Consequently, we suppose that $k\geq 2$.
We will focus on proving that $x_{n-2} = k-1$, since the remaining cases can be demonstrated using a similar argument. 
Note that $x$ is a revised ascent sequence, that is, $x\in \Cay_n$ with 
$\Ascbot(x) = \nub(x)$.
Since $x_{n-1} = k<m=x_n$, then $n-1 \in \Ascbot(x) = \nub(x)$.
This means that $x_{n-1}$ is the leftmost copy of $k$ in $x$,  it is also the 
unique copy of $k$ in $x$.
Suppose that the leftmost $k-1$ appears at the position $j$ in $x$.
It follows that $x_j = k-1 < x_{j+1}$.
If $j \leq n-3$, then $x_{j+1} >k$ as there exists only one copy of $k$ in $x$.
This yields that $x_jx_{j+1}x_{n-1}$ forms an occurrence of the pattern $132$ in $x$,  a contradiction to the fact that $x$ is $132$-avoiding.
Thus, $j = n-2$, namely, $x_{n-2} = k-1$.
This proves the claim.
That is, $x$ has the form $x = y12\cdots km$.
It is evident  that $\std(y)$ is a revised ascent sequence in $\RA_{n-1-k}(132)$.
Moreover, one can recover $x$ from $\std(y)$ easily.
Therefore, equation (\ref{equ:gs}) follows.
\qed

\begin{proposition}\label{prop:132}
	We have 
	\begin{equation}\label{equ:132}
		\RA_{132}(t) = \sum_{n\geq 1}|\RA_n(132)|t^n = \frac{t^{2} - 2 \, t+1 - \sqrt{t^{4} + 2 \, t^{2} - 4 \, t + 1} }{2 \, t}.
	\end{equation}
\end{proposition}

\pf
Let $G(t)$, $R(t)$, and $S(t)$ be the generating functions of $g_n$, $r_n$, and 
$s_n$, respectively.
%Turning the recurrence relations (\ref{equ:rs}), (\ref{equ:gr}), and (\ref{equ:gs}) into generating functions 
Multiplying both sides of equations (\ref{equ:rs}), (\ref{equ:gr}), and (\ref{equ:gs}) by $t^n$ and summing over all possible $n$ yields 
\[
\left\{
\begin{aligned}
	R(t) &= t+ tR(t) + S(t), \\
	G(t) &= \frac{t}{1-t} + \frac{G(t)}{t}\big(R(t) - \frac{t}{1-t}\big),\\
	S(t) &=  \frac{t^2G(t)}{1-t}.
\end{aligned}
\right.
\]
Note that $\RA_{132}(t) = G(t)$.
Solving the above system of equations for $G(t)$ gives (\ref{equ:132}).
\qed

Solving the above system of equations for $S(t)$ gives 
\[S(t)  = \frac{{\left(-t^{2} + 2 \, t - 1 + \sqrt{t^{4} + 2 \, t^{2} - 4 \, t + 1} \right)} \,t}{2 \, {\left(t- 1\right)}}.\]
Comparing it with the expression  of $\RA_{123}(t)$ as shown in (\ref{equ:123}), 
we derive that $$\RA_{123}(t) = t + \frac{S(t)}{t}.$$
This implies that $|\RA_n(123)| = s_{n+1}$ for $n\geq 2$.

Equation \eqref{equ:132} gives the sequence 
\[
(|\RA_{n}(132)|)_{n\geq 1} = (1,1,2,5,13,35,97,275,794, 2327,\cdots),
\]
which appears as A082582 in the OEIS \cite{oeis}.
This sequence  counts 
the Dyck paths of semilength $n$ with no {\bf uudd}.
%This sequence  also counts primitive modified ascent sequences that avoid
%the pattern $1243$ (or avoid the pattern $1123$) \cite{Zhou}.

Combining Lemma \ref{lem:main} and Proposition \ref{prop:132}, we have 
the following conclusion.

\begin{proposition}
		We have 
	\begin{equation*}
		\RA_{3132}(t) = \sum_{n\geq 1}|\RA_n(3132)|t^n = \frac{t^{2} - 2 \, t+1 - \sqrt{t^{4} + 2 \, t^{2} - 4 \, t + 1} }{2 \, t}.
	\end{equation*}
\end{proposition}

\section{Final remarks and future directions}\label{sec:final}

In this paper, we introduced the concept of revised ascent sequences and systematically investigated the enumeration of these sequences avoiding a 
single pattern.
We have derived enumeration results for revised ascent sequences that avoid all patterns of length $3$ except for the pattern $111$. 
The enumeration concerning the pattern $111$ remains open 
on both ascent sequences and modified ascent sequences too.
In addition to this, we also outline the following directions and suggestions for further research.

In Section \ref{sec:bijection}, we constructed a bijection between  ascent sequences and revised ascent sequences.
Extending this work, it would be particularly interesting to explore direct bijections between revised ascent sequences and other Fishburn structures, along with the corresponding statistics among these structures.

In Section \ref{subsection:123}, we proved that $\RA_n(123)$ are 
equinumerous with the set of Dyck paths of semilength $n$ with no {\bf uudu} (or with no {\bf uduu}).
Similarly, in Section \ref{subsection:132}, we proved that $\RA_n(132)$ are 
equinumerous with the set of Dyck paths of semilength $n$ with no {\bf uudd}.
One may ask whether there exist bijections among these structures. Additionally, it would be fascinating to investigate if there is a specific type of restricted Dyck paths that have a one-to-one correspondence with revised ascent sequences. 
Systematically studying pattern-avoiding revised ascent sequences and their corresponding pattern-avoiding Dyck paths could provide valuable insights and would be a particularly compelling area of research.

Given two patterns $\sigma$ and $\tau$, they are said to be {\em Wilf-equivalent}
if $|\RA_n(\sigma)| = |\RA_n(\tau)|$ hold for any $n\geq 1$.
Numerical evidences shows that there may exist a large number of Wilf-equivalences on revised ascent sequences.
Can we find more general conclusions about Wilf-equivalences on revised ascent sequences, similar to those presented in Lemma \ref{lem:main}.
Given an endofunction $x = x_1x_2\cdots x_n$, we say $x$ is 
\emph{primitive} if it contains no flat steps, where a {\em flat} step means  two consecutive equal entries $x_i = x_{i+1}$.
It may also be a worthwhile research direction to study the primitive revised 
ascent sequences and pattern avoidance in them.

It has been observed that revised ascent sequences avoiding patterns of length $4$ or longer also exhibit rich enumerative results and have close connections with numerous other combinatorial structures.
The detailed findings and conclusions related to this topic
are presently under preparation.
% will be presented in a separate paper.

	\section*{Acknowledgments}
	%The authors are very grateful to the referee for valuable comments and suggestions
	%which helped to improve the presentation of the paper.
	The work  was supported by
	the National Natural Science Foundation of China (11801378).
	
%	\section*{Data availability statements}
%	Data sharing is not applicable to this article as no datasets were generated or analyzed during the current study.


\begin{thebibliography}{100}
		
	\bibitem{Benyi}	
	B. B\'enyi, A. Claesson, and M. Dukes, 
	Weak ascent sequences and related combinatorial structures,
	{\em European J. Combin.}, {\bf 108} (2023), No. 103633.
	
		\bibitem{Benyi2024}
	B. B\'enyi, T. Mansour, and J.L. Ram\'irez, Pattern avoidance in weak ascent sequences,
	\emph{Discret. Math. Theor. Comput. Sci.}, {\bf 26} (2024).
	
	\bibitem{Bousquet}	
	M. Bousquet-M\'elou, A. Claesson, M. Dukes and  S. Kitaev,
	$(2+2)$-free posets, ascent sequences and pattern avoiding permutations,
	{\em J. Combin. Theory Ser. A}, {\bf 117} (2010), 884--909.
	
	\bibitem{Cerbai23a}
	G. Cerbai and A. Claesson, Fishburn trees, {\em Adv. Appl. Math.}, {\bf 151} (2023), Article 102592.
	
	\bibitem{Cerbai23b}
	G. Cerbai and A. Claesson, Transport of patterns by Burge transpose, 	{\em European J. Combin.}, {\bf 108} (2023), Article  103630.
	
	

	
	\bibitem{Cerb}
	G. Cerbai, Modified ascent sequences and bell numbers, 
	{\em Electron. J. Combin.}, {\bf 31(4)} (2024).
%	\href{https://arxiv.org/abs/2305.10820}{arXiv:2305.10820v1}, 2023.
	
	\bibitem{Cera}
	G. Cerbai, Pattern-avoiding modified ascent sequences, 
	\href{https://arxiv.org/abs/2401.10027}{arXiv:2401.10027v1}, 2024.
	
	\bibitem{Cerb2406}
	G. Cerbai, A. Claesson, and B.E. Sagan, 
	Modified difference ascent sequences and Fishburn structures,
	\href{https://arxiv.org/abs/2406.12610}{arXiv.org/abs/2406.12610}, 2024.
	
	\bibitem{Cerb2408}
	G. Cerbai, A. Claesson, and B.E. Sagan, 
	Self-modified difference ascent sequences,
	\href{https://arxiv.org/abs/2408.06959}{arXiv.org/abs/2408.06959}, 2024.
	
	\bibitem{Chen2013}
	W.Y.C. Chen, A.Y.L, Dai, T. Dokos, T. Dwyer, and B.E. Sagan,
	On $021$-avoiding ascent sequences,
	{\em Electron. J. Combin.}, {\bf 29(4)} (2013), \#P4.25.
	
	\bibitem{Yan2019}
	D.D. Chen, S.H.F. Yan, and R.D.P. Zhou, 
	Equidistribution statistics on Fishburn matrices and permutations,
	{\em Electron. J. Combin.}, {\bf 26} (2019), \#P1.11.
	
%	\bibitem{Claesson1}
%	A. Claesson, M. Dukes, and S. Kitaev, 
%	A direct encoding of Stoimenow's matchings as ascent sequences,
%	{\em Australas. J. Combin.}, {\bf 49} (2011), 47--59.
	
	\bibitem{Claesson2}
	A. Claesson and S. Linusson,
	$n!$ matching, $n!$ posets,
	{\em Proc. Amer. Math. Soc.}, {\bf 139} (2011), 435--449.
	
%	\bibitem{Conway}
%	A.R. Conway, M. Conway, A.E. Price, and A.J. Guttmann, 
%	Pattern-avoiding ascent sequences of length $3$, 
%	{\em Electron. J. Combin.}, {\bf 20(1)} (2022), \#P4.25

	\bibitem{Dukes2010}
	M. Dukes and R. Parviainen, Ascent sequences and upper triangular matrices
	containing non-negative integers,
	{\em Electron.   J.	Combin.},  {\bf 17} (2010),  R53.
	
%	\bibitem{Dukes2016}
%	M. Dukes, Generalized ballot sequences are ascent sequences,
%	{\em Australas. J. Combin.}, {\bf 64} (2016), 61--63.
	
	\bibitem{Dukes2019}
	M. Dukes and P.R.W. McNamara, 
	Refining the bijections among ascent sequences, (2+2)-free posets, integer matrices and pattern-avoiding permutations, 
	{\em J. Combin. Theory Ser. A}, {\bf 167} (2019), 403--430.
	
	\bibitem{Dukes2023}
	M. Dukes and B.E. Sagan, 
	Difference ascent sequences, 
	{\em Adv. Appl. Math.}, {\bf 159} (2024), Article 102736.
%	\href{https://arXiv:2311.15370v1}{arXiv:2311.15370v1}, 2023.
	
	
	
	
	\bibitem{Duncan}
	P. Duncan and E. Steingr\'imsson, 
	Pattern avoidance in ascent sequences,
	{\em Electron. J. Combin.}, {\bf 18} (2011), Article 226.
	
%	\bibitem{Fishburn1}
%	P.C. Fishburn, Intransitive indifference with unequal indifference intervals,
%	{\em J. Math. Psych.},  {\bf 7} (1970), 144--149.
%	
%	\bibitem{Fishburn2}
%	P.C. Fishburn, Inteval graphs and interval orders,
%	{\em Discrete Math.},  {\bf 55} (1985), 135--149.
%	
%	\bibitem{Fishburn3}
%	P.C. Fishburn, Interval orders and interval graphs:  a study of
%	partially ordered sets,
%	{\em John Wiley $\&$ Sons}, 1985.
	
	\bibitem{Levande}
	P. Levande, Fishburn diagrams, Fishburn numbers and their refined generating functions, {\em J. Combin.	Theory Ser. A},  {\bf 120} (2013),  194--217.
	
	\bibitem{Lin} Z. Lin, Patterns of relation triples in inversion and ascent sequences, {\em Theor. Comput. Sci.}, {\bf804} (2020), 115--125.
	
\bibitem{LY} Z. Lin and S.H.F. Yan, Vincular patterns in inversion sequences, {\em Appl. Math. Comput.}, {\bf364} (2020), Article 124672.

\bibitem{oeis} 
OEIS Foundation Inc., The On-Line Encyclopedia of Integer Sequences,  \href{http://oeis.org}{http://oeis.org}, 2011.
	
%	\bibitem{Mansour}
%	T. Mansour and M. Shattuck, Some enumerative results related to ascent sequences, 
%	{\em Discrete Math.},  {\bf 315} (2014), 29--41.
	
%	\bibitem{MS}
%	M. Martinez and C.D. Savage, Patterns in inversion sequences II: Inversion sequences avoiding triples of relations, {\em J. Integer Seq.}, {\bf21} (2018), Article 18.2.2.
	
	\bibitem{Stanley}
	R.P. Stanley,  
	Enumerative Combinatorics, vol. 1, second edition,
	Cambridge University Press, Cambridge, 2012.

	\bibitem{Stoimenow}
	A. Stoimenow, Enumeration of chord diagrams and an upper bound for Vassiliev invariants,
	{\em J. Knot Theory Ramifications}, {\bf 7} (1998), 93--114.

	
	
	\bibitem{West}
	J. West, Generating trees and forbidden subsequences, \emph{Discrete Math.}, {\bf 157} (1996), 363--374.
	
%	\bibitem{Zagier}
%	D. Zagier, Vassiliev invariants and a strange identity related to the dedekind eta-function, {\em Topology},
%	{\bf 40} (2001), 945--960.

   \bibitem{Yan}
   S.H.F. Yan, Ascent sequences and $3$-nonnesting set partitions, \emph{European J. Combin.}, {\bf 39} (2014), 80--94.
	
	\bibitem{Zhou}
	Y.C. Zang and R.D.P. Zhou,  Difference ascent sequences and related combinatorial structures, \href{https://arxiv.org/abs/2405.03275}{arXiv.org/abs/2405.03275}, 2024.
	
	

		
	\end{thebibliography}
\end{document}